\documentclass[12pt]{article}

\usepackage{amsmath}
\usepackage{amsthm}
\usepackage{amsfonts}
\usepackage{url}
\usepackage{xcolor}
\usepackage[backref=page]{hyperref} 
\hypersetup{
    colorlinks,
    linkcolor={red!50!black},
    citecolor={blue!50!black},
    urlcolor={blue!80!black},
}
\usepackage[lined,boxed,linesnumbered,algoruled]{algorithm2e}
\usepackage{cleveref}
\usepackage{authblk}
\usepackage[margin=1in]{geometry}
\usepackage{multirow}
\usepackage{float} 
\usepackage[all]{nowidow}
\usepackage{tikz-cd}

\newcommand{\mr}{\mathbb{R}}

\newcommand{\mc}{\mathbb{C}}

\newcommand{\bo}[1]{\mathbf{#1}}
\newcommand{\V}{\mathcal{V}}

\DeclareMathOperator{\rank}{rank}
\DeclareMathOperator{\im}{im}

\theoremstyle{definition}
\newtheorem{theorem}{Theorem}[section]
\newtheorem{corollary}[theorem]{Corollary}
\newtheorem{proposition}[theorem]{Proposition}
\newtheorem{definition}[theorem]{Definition}

\newtheorem{lemma}[theorem]{Lemma}
\newtheorem{example}[theorem]{Example}
\newtheorem*{example*}{Example}
\theoremstyle{remark}
\newtheorem{remark}[theorem]{Remark}

\theoremstyle{definition}
\newtheorem*{definition*}{Definition}
\theoremstyle{definition}
\newtheorem*{theorem*}{Theorem}
\theoremstyle{plain}
\newtheorem*{corollary*}{Corollary}
\theoremstyle{plain}
\newtheorem*{proposition*}{Proposition}
\theoremstyle{plain}
\newtheorem*{conj*}{Conjecture}
\newcommand{\bC}{\mathbb{C}}
\newcommand{\bR}{\mathbb{R}}
\newcommand{\sL}{\mathcal{L}}
\newcommand{\sU}{\mathcal{U}}
\newcommand{\starx}{{x^*}}

\crefname{appendix}{Appendix}{Appendices}

\DeclareMathOperator{\act}{\mathbf{Act}}
\DeclareMathOperator{\group}{\mathbf{Group}}
\DeclareMathOperator{\loc}{Mon}
\DeclareMathOperator{\iso}{iso}

\newcommand{\wt}[1]{{\widetilde{#1}}}

\renewcommand{\setminus}{-}

\title{On computing local monodromy and the \\
numerical local irreducible decomposition}
\author[1]{Parker B. Edwards}
\author[2]{Jonathan D. Hauenstein}
\affil[1]{Department of Mathematics and Statistics\\ Florida Atlantic University, Boca Raton, FL, USA}
\affil[2]{Department of Applied and Computational Mathematics and Statistics\\ University of Notre Dame, Notre Dame, IN, USA}
\date{}

\begin{document}
\maketitle
\begin{abstract}\noindent  
Similarly to the global case, the local structure of a holomorphic subvariety 
at a given point is described by its local irreducible decomposition.  
Following the paradigm of numerical algebraic geometry, 
an algebraic subvariety at a point is 
represented by a numerical local irreducible decomposition comprised of a local witness set for each local irreducible component.  
The key requirement for obtaining a numerical local irreducible
decomposition is to compute the local monodromy action of a 
generic linear projection at the given point,
which is always well-defined on any small enough neighborhood.
We characterize some of the behavior of local monodromy action
of linear projection maps under analytic continuation,
allowing computations to be performed beyond a local neighborhood.
With this characterization, we present an algorithm to 
compute the local monodromy action and corresponding numerical local irreducible
decomposition for algebraic varieties. 
The results are illustrated using several examples 
facilitated by an implementation in an open source software package.  
\end{abstract}
\section{Introduction}
Theories to understand the geometry and topology of a space at its singular points comprise a major and ongoing area of mathematical study. 
Some classically studied aspects of singularity theory include local invariants, 
local monodromy groups, and integration over singular spaces
(for a general overview, see the book~\cite{SingularityTheory}).
Computational methods directed towards identifying singularities,
understanding local structure, and 
stratifying spaces increasingly arise in applications as well, e.g.,~\cite{bendich2012strat,LandauSingularities,EllipticPDEs,Helmer2023whitney,lim2023hades}. 
We will focus here on a classical setting: given a system of algebraic or analytic equations and a point which satisfies them as input, compute information about the local structure of the solution set at that point in $\bC^N$. 

Understanding the \emph{global} structure of a solution set of a system
of algebraic equations is the foundational problem of algebraic geometry.
Globally, the solution set can be decomposed into finitely many 
irreducible components.  Following the numerical algebraic geometric paradigm
(for a general overview, see the books~\cite{BertiniBook,SommeseWamplerBook}),
each irreducible component is represented by a witness set 
yielding a corresponding numerical irreducible decomposition.
One key property is that global irreducibility is maintained
under intersection by a general hyperplane for irreducible components of
dimension at least two.  Therefore, all computations associated
with deciding irreducibility of positive-dimensional components 
can be reduced down to the complex curve case. 

For a germ $\bo{V}$ of a complex algebraic subvariety, there is a natural
analog yielding the local irreducible decomposition.
In fact, the theory presented below applies to germs of holomorphic subvarieties, though we will consider examples with algebraic inputs. 
Holomorphic germs also have a unique local irreducible decomposition (e.g.,
see~\cite[Thm. II.B.7]{gunning1990vol2}).
An analog following the numerical algebraic geometric paradigm
via local witness sets and a numerical local irreducible decomposition
was described in~\cite{BHS16}.  Since that work did not consider
how to actually compute such a numerical local irreducible decomposition,
the following develops new theoretical results yielding
an algorithm with a theoretical guarantee of computing the 
correct numerical local irreducible decomposition.
A software package implementing the algorithm is available at 
\url{https://github.com/P-Edwards/LocalMonodromy.jl}.

\begin{figure}[!t]
    \centering
    \includegraphics[width=0.22\linewidth]{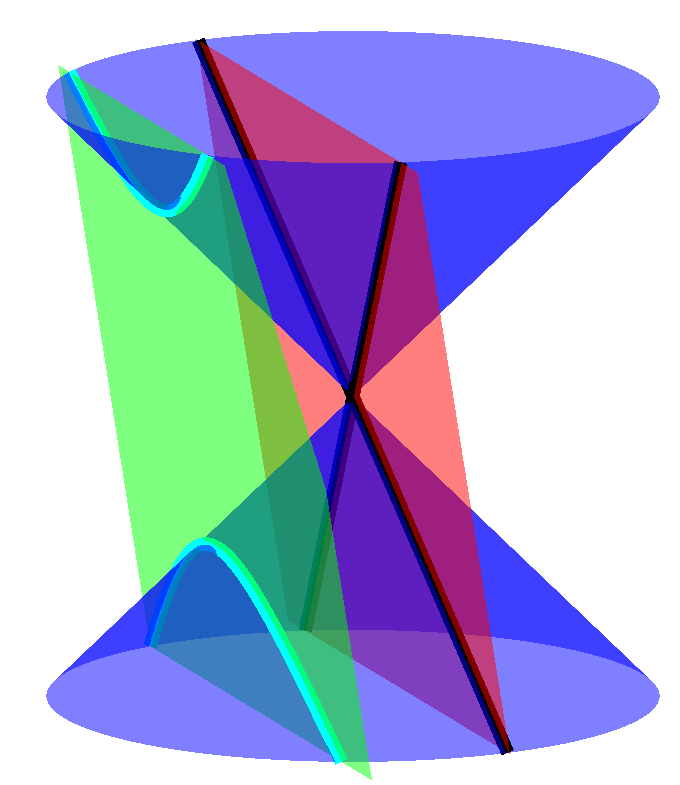} 
    \caption{Intersecting the cone (blue) with a general slice (green) yields an irreducible curve (cyan), while intersecting with a general slice through the origin (red) yields two lines (black).
    }
    \label{fig:Cone}
\end{figure}

One challenge with the local case is that one can not always reduce down to curves as in the global case.  To illustrate, consider the cone defined by $x_1^2+x_2^2-x_3^2=0$,
which is an irreducible surface in $\bC^3$.
Thus, for a general $\alpha\in\bC^3$, 
one can intersect this surface with a general hyperplane
defined by $x_1 = \alpha_1 x_2 + \alpha_2 x_3 + \alpha_3$ which yields an irreducible
curve, i.e., $(\alpha_1 x_2 + \alpha_2 x_3 + \alpha_3)^2 + x_2^2 - x_3^2$
is an irreducible polynomial for general $\alpha\in\bC^3$
as illustrated in Figure~\ref{fig:Cone}.
The origin is a singular point of the cone
in which the cone is locally irreducible at the origin.
However, for a general hyperplane of the form $x_1=\alpha_1x_2+\alpha_2x_3$ passing through the origin,
$(\alpha_1 x_2 + \alpha_2 x_3)^2 + x_2^2 - x_3^2$ is no longer irreducible
since every singular quadratic plane curve is 
simply a pair of intersecting lines as shown in Figure~\ref{fig:Cone}.

The following summarizes our theoretical results. 

\begin{theorem*}
If $\bo{V}$ is a reduced germ of a holomorphic subvariety of $\mc^N$ with pure dimension~$d$, there is a Zariski open set of linear projections $\mc^N\to\mc^d$ where, if $\wt{\pi}$ is a member, the projection map germ
$\boldsymbol{\wt\pi}\vert_\bo{V}:\bo{V}\to\bo{C^d}$ has a well-defined local monodromy action.  Moreover, if $\wt\pi\vert_V:V\to\mc^d$ is a proper projective holomorphic map with pure $(d-1)$-dimensional critical locus representing $\boldsymbol{\wt\pi}\vert_{\bo{V}}$, then the local monodromy action is a 
sub-action of the monodromy action of $\wt\pi\vert_V$ which can 
be computed using numerical algebraic geometry.
The corresponding local monodromy group decomposes into
orbits with one for each local irreducible component. 
\end{theorem*}

We show this in two parts. First, we apply a theorem of Hamm and {L\^{e} D\~{u}ng Tr\'{a}ng}~\cite{HL73} on the fundamental group of an analytic hypersurface to define a local monodromy action generated by loops contained in a complex line for a projection $\wt\pi$. This is a localized version of a strategy for global monodromy actions taken in~\cite{HRS18} using a theorem of Zariski~\cite{zariski1937}.
Although the global case can always be reduced to curves, this result permits a reduction
down to surfaces in the local case.  
Since the theorem of~\cite{HL73} is local, 
it requires finding a small enough restriction to localize
the monodromy computations.  We overcome this by 
characterizing how this local monodromy action includes in a simple way into the monodromy action of any appropriate analytic continuation 
of $\wt\pi$ thereby extending the domain and facilitating computations. 
The following example gives an overview of the approach. 

\begin{figure}[!t]
    \centering
    \includegraphics[width=0.29\linewidth]{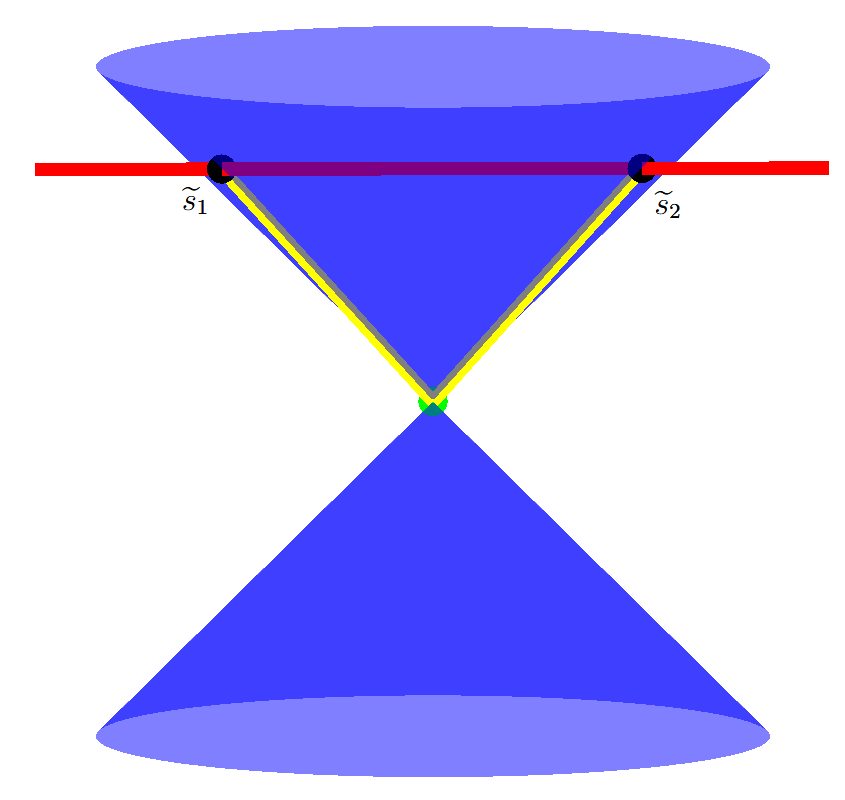}
    \caption{Intersection of the cone (blue) with line (red)
    yields two points (black) which are the start points of 
    two paths (yellow) that limit to the origin (green).
    }
    \label{fig:ConeWitness}
\end{figure}

\vspace{-0.1in}
\paragraph{Illustrating example}
Reconsider the cone $C\subset\bC^3$ defined by $f(x)=x_1^2+x_2^2-x_3^2=0$
at the origin which yields a reduced germ of $\bC^3$ with
pure dimension $2$.  For illustrative purposes, 
consider the sufficiently general
linear projection $\wt \pi(x) = (x_1+x_2,x_3)$.
Thus, for the sufficiently general point 
$\gamma_1 = (1/2,1)\in\bC^2$, 
the fiber
$\wt\pi^{-1}(\gamma_1)\cap C$ consists of two points.
Along the segment $\gamma(t) = t\gamma_1$ for $t\in(0,1]$,
$\wt\pi^{-1}(\gamma(t))\cap C$ defines two solution paths
starting at these two points.  In this case, both points ``localize'' by limiting to the origin as illustrated in Figure~\ref{fig:ConeWitness}. Call them $\wt s_1$ and $\wt s_2$. 

\begin{figure}[!t]
    \centering
    \begin{tabular}{ccc}
    \includegraphics[width=0.21\linewidth]{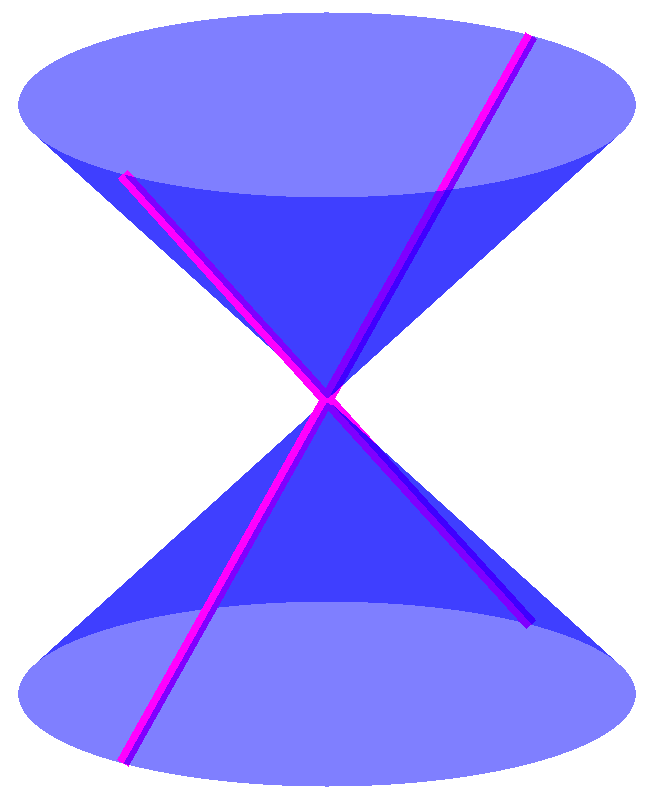} 
    & \hspace{1in} &  
    \includegraphics[width=0.25\linewidth]{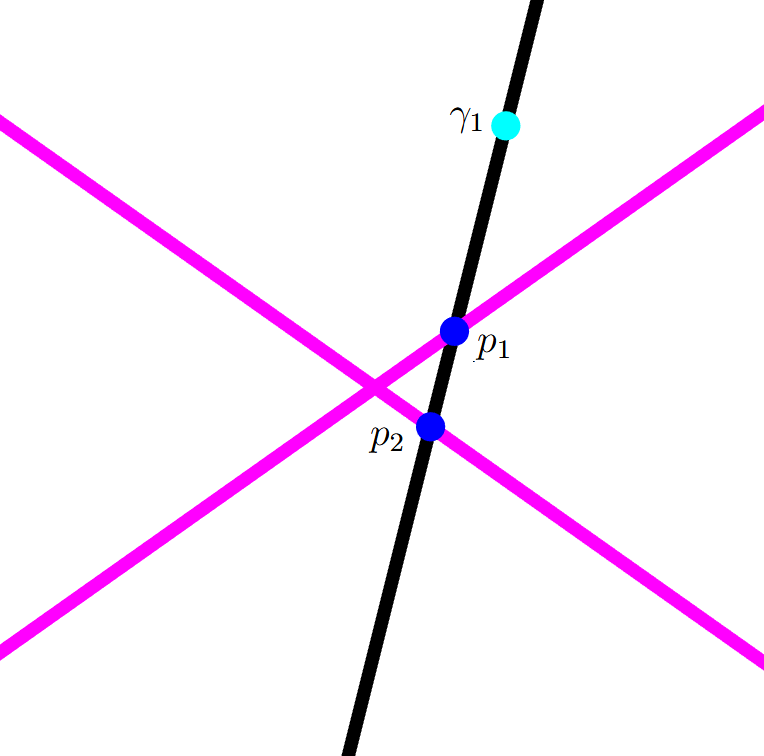}
    \\
    (a) && (b)  \\
    \end{tabular}
    \caption{(a) Critical points (magenta) with respect to $\wt \pi$
    on the cone (blue); (b) critical locus (magenta) 
    in the image of $\wt \pi$ intersected with a line (black)
    passing through $\gamma_1$ (cyan) yielding
    two points (blue).
    }
    \label{fig:ConeWitness2}
\end{figure}

Next, we need to consider the critical points
of $\wt \pi$, which comprise a hypersurface on the cone
as illustrated in Figure~\ref{fig:ConeWitness2}(a).
To compute the 
critical locus in terms of the projection~$\wt \pi$,
we consider another sufficiently general linear projection
$\pi:\bC^2\rightarrow\bC$ with $z\mapsto z_1-z_2/4$.
Let $\theta(1) = \pi(\gamma_1)$ and consider 
the linear space $\sL_{\theta(1)}\subset\bC^2$
defined by $\pi(z)=\theta(1)$.  Since the determinant
of the Jacobian matrix of $f$ and $\wt \pi$ is $2x_1-2x_2$,
we need to solve the system
\begin{equation}\label{eq:ConeCritical}
\hbox{\small $
\left[\begin{array}{c} f(x) \\ \hline \wt\pi(x) - z \\ \hline \det J(f,\wt\pi)(x) \\ \hline \pi(z) - \theta(1) \end{array}\right] = 
\left[\begin{array}{c}
x_1^2+x_2^2-x_3^2 \\ \hline x_1+x_2 - z_1 \\ x_3 - z_2 \\ \hline 2x_1-2x_2 \\ \hline z_1-z_2/4-1/4
\end{array}\right]
=0.$}
\end{equation}
In terms of the $z$ coordinates, this yields two points
as illustrated in Figure~\ref{fig:ConeWitness2}(b).
Replacing $\theta(1)$ with $\theta(t) = \pi(\gamma(t)) = t/6$
in \eqref{eq:ConeCritical}, yields two solution paths
starting at these two points.  Since, in this case,
both limit to the origin, we call the two points $p_1$ and $p_2$. 

\begin{figure}[!t]
    \centering
    \begin{tabular}{ccc}
    \includegraphics[width=0.2\linewidth]{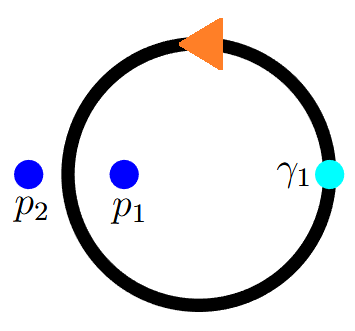} 
    & \hspace{1in} & 
    \includegraphics[width=0.22\linewidth]{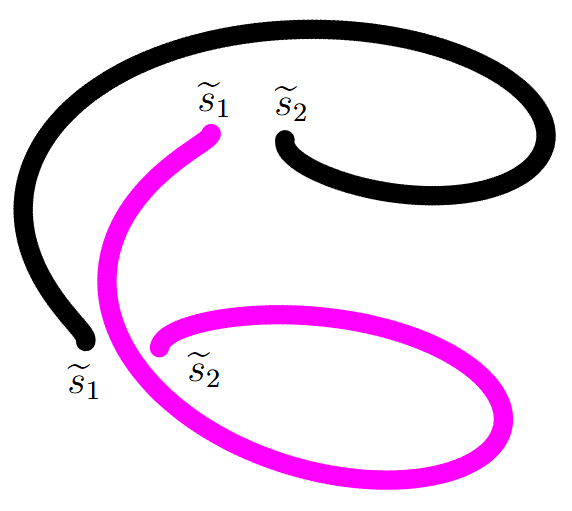}
    \\
    (a) && (b)  \\
    \end{tabular}
    \caption{(a) Illustration of a basic loop 
    starting at $\gamma_1$ and only encircling $p_1$ once counterclockwise; 
    (b)~pictorial illustration of monodromy action
    interchanging $\wt s_1$ and $\wt s_2$.
    }
    \label{fig:ConeLoop}
\end{figure}

Our theoretical results show that one can compute
the local monodromy action of the germ of the cone at the origin
by considering the local monodromy action 
arising from loops in $\sL_{\theta_1}\setminus\{p_1,p_2\}$
lifted to the cone.  
In particular, one can view $\sL_{\theta_1}\setminus\{p_1,p_2\}$
as $\bR^2$ with two points removed and observe that the action of 
a basic loop starting at $\gamma_1$ that only 
encircles~$p_1$ once counterclockwise, 
as illustrated in Figure~\ref{fig:ConeLoop}(a), suffices to generate the local monodromy action.  
Such a loop
lifts to two paths in $C$ starting at $\wt s_1$ and $\wt s_2$.
The corresponding path starting at $\wt s_1$ ends at $\wt s_2$
and vice versa as pictorially illustrated in Figure~\ref{fig:ConeLoop}(b).
Note that a basic loop starting at $\gamma_1$ that only 
encircles~$p_2$ once counterclockwise performs the same monodromy action as a basic loop encircling $p_1$ clockwise.
Hence, the local monodromy action
has a single orbit and the local monodromy group is the symmetric group on two elements,
which shows that the germ of the cone at the origin
is irreducible.

\medskip

The rest of the paper justifies this process 
and is organized as follows.
\Cref{sec:background} recalls some essential results and 
definitions about germs of holomorphic subvarieties and their local irreducible decompositions. \Cref{sec:theory} lays out our theory of local monodromy actions and groups for projection maps, which is used in \Cref{sec:algorithm} to 
justify an algorithm for computing these objects. 
\Cref{sec:examples} contains several examples computed using an implementation of the algorithm.  A short conclusion is provided in \Cref{sec:conclusion}.

\section{Background}\label{sec:background}
We review local parameterization of holomorphic subvarieties, give an overview of homotopy continuation and numerical algebraic geometry, and state Hamm and {L\^{e} D\~{u}ng Tr\'{a}ng}'s theorem on fundamental groups of hypersurface complements~\cite{HL73}. 

\subsection{Local irreducible decomposition and local parameterization}
A germ $\bo{V}$ of a holomorphic subvariety of an open set $\wt U\subseteq\mc^N$ at $x^*\in \wt U$ has an irreducible decomposition.
More precisely, using Gunning's notation~\cite[Thm. II.B.7]{gunning1990vol2}, 
$\bo{V}$ can be written as a finite union of germs $\bo{V} = \bo{V}_1 \cup \bo{V}_2 \cup \dots \cup \bo{V}_m$ where each $\bo{V}_j$ is an irreducible germ, $\bo{V}_j\not= \bo{V}$ if $m > 1$, 
$\bo{V}_j \not\subseteq\bo{V}_k$ for $j\not= k$, and the germs $\bo{V}_j$ are uniquely determined up to relabeling. 
Following \cite{BHS16}, the local irreducible decomposition of an irreducible subvariety (algebraic or holomorphic) $V$ at $x^*\in V$ is the decomposition given by the germ of $V$ at $x^*$. 
Moreover, if $V$ is reducible, then the local irreducible decomposition of 
$V$ at $x^*$ is the union of local irreducible decompositions of its (global) irreducible components. 

Since one may always translate $x^*$ to the origin, 
it suffices to consider germs at the origin, which we will do going forward to simplify notation except where otherwise indicated. 
While holomorphic subvarieties may exhibit more complicated global behavior, they exhibit the structure of finite branched coverings locally. The following collects several standard results making this precise in a useful format for our purposes. 

\begin{lemma}\label{lem:local}
    Let $\bo{V}$ be a pure $d$-dimensional germ of a holomorphic subvariety of $\mc^N$ with~$V$ a representative. There exists an (algebraically) Zariski open set of linear projections \mbox{$\mc^N\to\mc^d$} where the following holds 
    provided that $\pi:\mc^N\to\mc^d$ is a member. For all small enough open balls $\wt{B}\subseteq \mc^N$ and $B\subseteq\mc^d$ at the origin, $\hat{V}:= V\cap\wt{B}\cap\pi^{-1}(B)$ has the form $\hat{V} = \hat{V}_1 \cup \dots \cup \hat{V}_m$ where the $\hat{V}_i$ are irreducible holomorphic subvarieties representing the irreducible components of $\bo{V}$. Furthermore: 
    \begin{enumerate}
        \item $\pi\vert_{\hat{V}}$ and $\pi\vert_{\hat{V}_i}$ are finite branched holomorphic coverings of $B$ with $0\in\mc^N$ the only element of the fiber over $0\in\mc^d$ for both maps.
        \item The image of the branch locus for $\pi\vert_{\hat{V}}$ and $\pi\vert_{\hat{V}_i}$ is a holomorphic subvariety of $B$ with the same dimension as the branch locus.
        \item If $\wt{R}$ and $\wt{R}_i$ are the branch locuses of $\pi\vert_{\hat{V}}$ and $\pi\vert_{\hat{V}_i}$ respectively, the monodromy action on any fiber of $\pi\vert_{\hat{V}_i\setminus \wt{R}_i}$ is transitive and the monodromy action of $\pi\vert_{\hat{V}\setminus \wt{R}}$ on any fiber partitions the fiber into orbits, one for each local irreducible component. 
    \end{enumerate}
\end{lemma}
\begin{proof}
    The first statement follows from the local parameterization theorem, e.g., 
    see~\cite[Lem.~II.E.12]
    {gunning1990vol2}. Noting that the branch locus and its image are holomorphic subvarities for a 
    finite holomorphic branched covering, e.g., see~\cite[Thm.~II.C.13,14]{gunning1990vol2}, the second 
    statement follows from the first using Remmert's proper mapping theorem,
    e.g., see~\cite[II.N.1]{gunning1990vol2}. 
    Since the germ of $\hat{V}_i$ at the origin is irreducible, possibly shrinking $B$ further, we have that $\hat{V}_i\setminus \wt{R}_i$ is path connected~\cite[II.E.13]{gunning1990vol2} and the corresponding monodromy action is transitive. Any point in the intersection of two distinct irreducible components of $\hat{V}$ is a branch point of~$\pi\vert_{\hat{V}}$, so, in particular, every point in a fiber over a regular point of $\pi\vert_{\hat{V}}$ is contained in one and only one irreducible component $\hat{V}_i$.
\end{proof}

\begin{remark}\label{rem:ReplaceBranchCritical}
    Taking the balls $\wt{B}$ and $B$ in the above lemma small enough, we may assume that $
    \hat{V}$ is defined as the zero set of a system $F$ of holomorphic functions on $\wt{U}\subseteq\mc^N$. If $
    \hat{V}$ is a reduced complete intersection with respect to such a system, the critical 
    points of $\pi\vert_{\hat{V}}$ 
    and~$\pi\vert_{\hat{V}_i}$ 
    are either empty or 
    holomorphic subvarieties of 
    dimension $d-1$ which contain the corresponding branch points. In this case, 
    \Cref{lem:local} remains true replacing branch loci with the corresponding critical loci.
  The critical loci is
  commonly used to facilitate numerical algebraic geometric computations involving the branch loci, 
  e.g., see~\cite{GeometricGenus,HRS18}.
\end{remark}

\subsection{Numerical algebraic geometry and homotopy continuation}

Numerical algebraic geometry (for a general overview, see the books~\cite{BertiniBook,SommeseWamplerBook})
represents an irreducible algebraic variety via a witness set.
Suppose that $f:\bC^N\rightarrow\bC^n$ is a polynomial system and
$\V(f) = \{x\in\bC^N~|~f(x) = 0\}$.
If $V\subset\V(f)$ is irreducible of dimension $d$, 
then, for a general codimension $d$ linear space $\sL$, 
the intersection $V\cap\sL$ is finite and the number of such points in
the intersection is equal to $\deg V$.
The set $\{F,\sL,V\cap\sL\}$ is called a \emph{witness set} for $V$
and $V\cap\sL$ is called a \emph{witness point set}. 
The following considers a local version~\cite{BHS16}.

\begin{definition}[From \cite{BHS16}]\label{def:local_witness_sets}
    Let $f:\mc^N\to\mc^n$ be a system of functions which are holomorphic in a neighborhood of $\starx\in\mc^N$ with $f(\starx) = 0$. Let $V\subseteq\mc^N$ be a local irreducible component of $\mathcal{V}(f)$ at $\starx$ of dimension $d$ and $\wt{L}_1,\dots,\wt{L}_d:\mc^N\to\mc$ be general linear polynomials such that $\wt{L}_i(\starx)=0$. For $u\in\mc^d$, 
    consider the linear space $\sL_u = V(\wt{L}_1-u_1,\dots,\wt{L}_d-u_d) \subset \mc^N$.
    A \emph{local witness set} for $V$ is the triple $\{f,\sL_{u^*},W\}$ defined in a 
    sufficiently small neighborhood~$U\subset\mc^d$ of the origin for general $u^*\in U$ where $W$ is the finite subset of points of $V\cap\sL_{u^*}$ which converge to $\starx$ as $t\rightarrow0$ along any path defined by $V\cap\sL_{u(t)}$ such that $u:[0,1]\to U$ with $u(0)=0$ and $u(1)=u^*$.
\end{definition}

\begin{definition}[\cite{BHS16}]
    A \emph{numerical local irreducible decomposition} of a holomorphic subvariety $V$ at $\starx\in V$ is a formal union of local witness sets, one for a representative of each irreducible component of the germ $\bo{V}$ at $\starx$. 
\end{definition}

The set $W$ in \Cref{def:local_witness_sets} is called a {\em local witness point set}. Some conditions which~\cite{BHS16} leaves implicit are required for 
\Cref{def:local_witness_sets} to make sense. It is sufficient for the projection map $\pi:\mc^N\to\mc^d$ defined by the linear forms $\wt{L}_1,\dots,\wt{L}_d$ to 
have $\pi\vert_{V\cap \pi^{-1}(U)}$ be an (unbranched) covering map. 
For example, if $V$ is algebraic, the Noether normalization theorem shows that this is true for a Zariski open set of linear forms after removing a branch locus. 

\Cref{def:local_witness_sets} is dependent upon computing start
points of paths which converge to $x^*$.  This is an example of
the use of homotopy continuation in numerical algebraic geometry.
In particular, suppose that $\pi:V\to Y$ is 
a finite branched holomorphic covering map and $\gamma:[0,1]\to Y$ 
with $\gamma\vert_{(0,1]}$ a smooth path
into the regular part of $Y$. Then, $\pi^{-1}(\gamma(1))$ is finite with, say, $k$ points and $\gamma$ lifts through $\pi$ to a set of $k$ paths $\wt{\gamma}_i:[0,1]\to V$, each having $\wt{\gamma}_i(1)$ a distinct point of $\pi^{-1}(\gamma(1))$, e.g., see~\cite[Thm.~3]{morgansommese1989}.
 These paths are smooth and disjoint in the sense that each restriction $\wt{\gamma_i}\vert_{(0,1]}$ is smooth and none of the images of the $\wt{\gamma_i}\vert_{(0,1]}$ intersect. 
 Each lift $\wt{\gamma_i}$ is a solution to an initial value problem with initial values given by $\pi^{-1}(\gamma(1))$ and numerical homotopy continuation methods are designed to track the solutions numerically. 
 The lifts $\wt{\gamma_i}$ are sometimes called \emph{solution paths}.

Using homotopy continuation to move linear slices in a witness
set is a powerful tool in numerical algebraic geometry.
In the current context, we can consider general loops 
in the corresponding Grassmannian to induce a monodromy action
that can be used to identify the global and local irreducible components
by partitioning the fiber into orbits, one for each global
and local irreducible component, respectively.

\begin{example}\label{ex:Cone}
Globally, the cone $V=\V(x_1^2+x_2^2-x_3^2)\subset\bC^3$
is irreducible of dimension $2$ and degree $2$.
Consider the sufficiently
general line \mbox{$L=\V(x_1+x_2-1/2, x_3-1)\subset\bC^3$}
as in the Introduction.
Figure~\ref{fig:ConeWitness} 
illustrates the witness point set $V\cap L$.
Since $V$ is irreducible, the monodromy group obtained by permuting the two
points $V\cap L$ along general loops in the Grassmannian of lines in $\bC^3$
starting and ending at $L$ is the symmetric group on~two~elements.

Let $L_t=\V(x_1+x_2-t/2, x_3-t)\subset\bC^3$
so that $V\cap L_t$ defines $2$ paths starting from the two points in 
$V\cap L$.  
As also shown in Figure~\ref{fig:ConeWitness}, both paths 
converge to the origin, which means that a local witness set
for $V$ at the origin will have a local witness point set consisting
of~$2$ points.  
As shown in the Introduction, the cone is irreducible
at the origin so that the corresponding local monodromy group is also the symmetric group on~two~elements.  
\end{example}

\subsection{Hyperplane sections and fundamental groups of complements}

As mentioned in the Introduction, Zariski's theorem~\cite{zariski1937}
was used in \cite{HRS18} to compute monodromy groups
via a surjection of fundamental groups under slicing.  
In order to consider the local case, let 
$\wt B_r$ denote the ball of radius $r$ in $\mc^N$ centered at the origin and $B_r$ denote the same in $\mc^d$ for $r>0$.
A ``Zariski theorem of Lefschetz type'' due to Hamm and {L\^{e} D\~{u}ng Tr\'{a}ng}~\cite{HL73} shows that the map induced by inclusion 
of fundamental groups 
$\pi_1((B_\rho\setminus R)\cap \ell,b)\to \pi_1(B_\rho\setminus R,b)$ is surjective for a holomorphic hypersurface $R$, sufficiently small $\rho$, and any well-behaved complex line $\ell\subseteq\mc^d$. For hyperplane $L\subseteq\mc^d$, by abuse of notation, let $L$ also denote some fixed linear form such that $L=0$ defines the hyperplane $L$.
Let $L^\theta=\V(L-\theta)$ for~$\theta\in\bC$.

\begin{lemma} Let $R$ be a pure $(d-1)$-dimensional holomorphic subvariety of $\mc^d$ containing the origin. There exists a Zariski open subset of hyperplanes $\sU\subseteq \text{Gr}(d-1,d)^{d-1}$ where, for all $(L_1,L_2,\dots,L_{d-1}) \in \sU$, there exists $A >0$ such that, for all $\rho$ with $0 < \rho \leq A$, there is $\theta(\rho) > 0$ such that 
if $0 < \vert\theta\vert \leq \theta(\rho)$, the homomorphism induced by inclusion 
\[\pi_1((B_\rho \setminus R)\cap_{i=1}^{d-2} L_i \cap L^\theta_{d-1},b) \to \pi_1(B_\rho\setminus R,b)\] is surjective for any $b\in (B_\rho\setminus R)\cap_{i=1}^{d-2} L_i \cap L^\theta_{d-1}$.
\label{lem:lefschetz}
\end{lemma}
\begin{proof}
This is simply repeated application of~\cite[Thm.~0.2.1]{HL73}. 
Each application requires that $R$ be defined as the vanishing locus of a single holomorphic function. This is true for small enough $\rho$ for pure $(d-1)$-dimensional $R$, e.g., see~\cite[II.G.5]{gunning1990vol2}.
\end{proof}
\begin{remark}
    If $d = 2$, the indicated intersection of hyperplanes $\cap_{i=1}^{d-2} L_i$ should be understood to denote all of $\mc^2$. 
    When $d=1$, the intersection $\cap_{i=1}^{d-2} L_i \cap L_{d-1}^\theta$ should be understood as $\mc$ in which case $B_\rho\setminus R$ is a disc with a single puncture at the origin for small enough $\rho$. 
\end{remark}
Continuing the same notation with $\ell=\cap_{i=1}^{d-2} L_i \cap L_{d-1}^\theta$ for sufficiently small $\theta$, note that $(B_\rho\setminus R)\cap\ell = (B_\rho\cap\ell)\setminus (B_\rho\cap R\cap \ell)$. For general $\ell$, possibly shrinking $B_\rho$ further, one has from \Cref{lem:local} that $B_\rho\cap R\cap\ell$ is finite. Subsequently, $(B_\rho\setminus R)\cap\ell$ is homeomorphic to a real plane with the finitely many points $B_\rho\cap R\cap\ell$ missing. Its fundamental group is generated by a set of homotopy classes of based loops, one per missing point. Each loop in such a set encircles exactly one point in $B_\rho\cap R\cap\ell$ once counterclockwise. 

\section{Local monodromy actions and groups}\label{sec:theory}

Lifting from loops that generate the fundamental group, 
\Cref{lem:lefschetz} immediately yields 
that the monodromy actions for representatives of the germ $\bo{V}$ stabilize in small enough neighborhoods. Denote any intersection of the form $\cap_{i=1}^{d-2} L_i \cap L^\theta_{d-1}$ as in that result by~$\sL_{\theta}$. We will require some machinery to define a limiting process that behaves well with respect to homotopy continuation methods. The following setup definition is motivated by the requirements of one of Morgan and Sommese's foundational parameter homotopy continuation theorems~\cite[Thm.~3]{morgansommese1989}. 

\begin{definition}\label{def:mon_data}
    A \emph{monodromy representative} for a germ of a holomorphic map $\boldsymbol{\wt\pi}:\bo{V}\to\bo{C}^d$ is comprised of the following sets, maps, and commutative diagram, where $\wt\pi$ is a representative of $\boldsymbol{\wt\pi}$:
    \[ 
    \begin{tikzcd}
    \mc^N \arrow[r,"\widetilde{\pi}"] & \mc^d \arrow[r,"\pi"] & \mc^{d-1} \\
    V \arrow[u,hook,"\subseteq"] \arrow[ur] & R \arrow[u,hook,"\subseteq"] \arrow[ur] &  \\
    \wt{R} \arrow[u,hook,"\subseteq"] \arrow[ru,two heads] & &
    \end{tikzcd}
    \]
    Furthermore:
    \begin{itemize}
        \item $V$ is a pure $d$-dimensional holomorphic subvariety of $\mc^N$ containing the origin.    
        \item $\wt{R}$ is a pure $(d-1)$-dimensional holomorphic subvariety of $V$ containing the origin. 
        \item $\wt{\pi}\vert_V$ is a proper projective holomorphic map and the restriction $\wt{\pi}\vert_{V\setminus \wt{R}}$ is a finite unbranched covering of its image. 
        \item $R = \wt{\pi}(\wt{R})$ is a pure $(d-1)$-dimensional holomorphic subvariety of $\mc^d$.
        \item The restriction $\pi\vert_R$ is a proper projective holomorphic map. 
    \end{itemize}
\end{definition}

\begin{remark} We will reuse the notation in \Cref{def:mon_data} going forward to refer to the components of any given monodromy representative.
\end{remark}

\begin{remark}\label{rem:key_example}
    Let $\bo{V}$ be a pure $d$-dimensional germ of a holomorphic subvariety at the origin in $\mc^N$. As 
    discussed in the previous section, for linear projection $\wt\pi:\mc^N\to\mc^d$ satisfying \Cref{lem:local}, there is a representative $\hat{V}$ of $\bo{V}$ where $\wt{\pi}\vert_{\hat{V}}$ is a finite branched holomorphic covering of a ball $B$ centered at $0\in\mc^d$ and $0\in\mc^N$ is the unique element in $\wt{\pi}^{-1}(0)$. If the origin is a singular point of $V$, the map $\wt{\pi}\vert_{\hat{V}}$ has a nonempty branch locus $\wt{R}\subseteq\hat{V}$ which 
    has codimension at least $1$ in $\hat{V}$. 
If the codimension is larger than $1$, 
simply connectedness yields trivial
local monodromy.  
Moreover, as in Remark~\ref{rem:ReplaceBranchCritical},
one can always replace the branch locus
with the critical locus which, by abuse of notation, we will also call $\wt{R}$.
Therefore, without loss of generality, we can assume that $\wt{R}$ has codimension $1$ in $\hat{V}$.
    
    Possibly shrinking $B$ further, let $L_1,\dots,L_{d-1} \subseteq\mc^d$ be hyperplanes for which \Cref{lem:lefschetz} applies to $R$ and $B$, and such that the projection $\pi:\mc^d\to\mc^{d-1}$ defined by $(L_1,L_2,\dots,L_{d-1})$ is a finite proper branched covering of its image. This is true for generic choices of hyperplanes since \Cref{lem:lefschetz} and the local parameterization theorem both apply generically. Then, the diagram in \Cref{def:mon_data} with $V$ replaced by $\hat{V}$ and other notation referring to the specific choices in this example is a monodromy representative for $\boldsymbol{\wt{\pi}\vert_\bo{V}}$.     
\end{remark}

\begin{definition}
    Given a pure $d$-dimensional germ $\bo{V}$ of holomorphic subvariety of $\mc^N$, call any monodromy representative of a linear projection $\boldsymbol{\wt{\pi}}:\bo{V}\to\bo{C}^d$ which fulfills the conditions in \Cref{rem:key_example} a \emph{localized monodromy representative} for $\bo{V}$.
\end{definition}

\begin{example} \label{ex:algebraic_continuation}
Let $V=\V(f)\subseteq\mc^N$ be a pure $d$-dimensional algebraic complete intersection where $f:\mc^N\to\mc^{N-d}$ is algebraic.
Let $Jf(x)$ denote the Jacobian matrix of $f$ at $x\in\bC^N$.
Assume that $V$ is reduced with respect to $f$ in the sense that $N-\rank(Jf(v))=d$ at every regular point $v\in V$ 
and suppose that $V$ has a singular point at the origin. 
Let $\wt{\pi}:\mc^N\to\mc^d$ be a linear projection 
where $\wt{\pi}\vert_V$ represents the germ of a local parameterization at the origin in the sense of \Cref{lem:local}. 
Hence, the critical locus of $\wt\pi$, which contains the branch locus, 
is either empty
or of pure dimension $d-1$ since 
it is characterized by satisfying the critical equations: $f(z) = 0$ and $\det J(f,\wt\pi)(z) = 0$.

Typically, it is challenging to compute equations defining the image of the critical locus under $\wt\pi$ directly. 
Instead, consider the graph $G := \{ (\wt{v},v)\in \mc^N\times\mc^d \mid f(\wt{v})=0,\wt{\pi}(\wt{v}) = v\}$. 
The critical points of $\wt\pi\vert_V$ correspond to the points on this graph defined by 
\[CG := \{(\wt{v},v)\in G \mid \rank J(F,\wt{\pi})(\wt{v}) < N \}.\] Define $\wt{R}$ to be the image of the projection of $CG$ onto the first factor with $R = \wt\pi(\wt{R})$ being 
the projection onto the second factor. 
Note that $\wt{R}$ is the critical locus of $\wt\pi\vert_V$
which contains the branch locus.
The elements of \Cref{def:mon_data} with notation referring to the specific choices in this example are a monodromy representative for $\boldsymbol{\wt\pi\vert_V}$.
\end{example}

Our limiting procedure for defining local monodromy groups in this context will proceed by tracking the monodromy action defined by a monodromy representative along a continuous path $\gamma:[0,a]\to\mc^d$ that goes from 0 to some other point. The constraints we need to place on $\gamma$ are determined again by parameter homotopy continuation considerations. For convenience, the following lemma collects those constraints and is a direct corollary of~\cite[Thm.~3-(2,4)]{morgansommese1989}. 
The solution paths arising can be tracked using homotopy continuation.
\begin{lemma}\label{lem:lifting}
    Given a monodromy representative as in \Cref{def:mon_data}, let $\gamma:[0,a]\to\mc^d$ be a continuous path with $\gamma\vert_{(0,a]}$ smooth, $\im(\gamma)\subseteq\wt{\pi}(V)$. There exist Zariski open dense subsets $V_0\subseteq \wt{\pi}(V)$ and $R_0\subseteq\pi(R)$ where: 
    \begin{enumerate}
        \item If $K_1$ is the line containing $\gamma(0)$ and $\gamma(a)$, then there are at most countably many points in $K_1\setminus(K_1\cap V_0)$ and they are geometrically isolated. If also $\gamma((0,a])\subseteq K_1\cap V_0$, then $\gamma$ lifts to $V$ through $\wt{\pi}\vert_V$ as a disjoint set of finitely many smooth paths $[0,a]\to V$.
        \item If additionally $\im(\pi\circ\gamma) \subseteq \pi(R)$ and $K_2$ is the line in $\mc^{d-1}$ containing $\pi(\gamma(0))$ and $\pi(\gamma(a))$, then there are at most countably many points in $K_2\setminus (K_2\cap R_0)$ and they are geometrically isolated. If also $(\pi\circ\gamma)((0,a]) \subseteq K_2\cap R_0$, then $\pi\circ\gamma$ lifts through $\pi\vert_{R}$ as a disjoint set of finitely many smooth paths $[0,a]\to R$.
    \end{enumerate}
\end{lemma}

\begin{definition}
For a monodromy representative as in \Cref{def:mon_data}, a path 
\mbox{$\gamma:[0,a]\to\mc^d$} starting at $0$ is a \emph{limiting path} for that representative if it fulfills all the conditions from 
    \Cref{lem:lifting}, including those in Items 1 and 2, if $\im(\gamma)\subseteq \cap_{i=1}^{d-2} L_i$, and $\Vert \gamma \Vert:[0,a]\to\mr$ is an increasing function.
\end{definition}

\begin{definition}
    Let $\wt{\epsilon},\rho > 0$. The \emph{restriction} of a monodromy representative with limiting path $\gamma:[0,a]\to\mc^d$ to $\wt{B}_\wt{\epsilon}$ and $B_\rho$ is obtained by replacing: 
    \begin{itemize}
        \item $V$ with $V\cap \wt{B}_\wt{\epsilon} \cap \wt{\pi}^{-1}(B_\rho)$
        \item $\wt{R}$ with $\wt{R}' = \wt{R}\cap \wt{B}_\wt{\epsilon}\cap\wt{\pi}^{-1}(B_\rho)$
        \item $R$ with $R' = \wt{\pi}(\wt{R}')$
        \item $\gamma$ with $\gamma\vert_{[0,\beta]}$ where $\beta = \sup( \{ t\in [0,a] \mid \gamma([0,t])\subseteq B_\rho\cap \pi^{-1}(\pi(R')) \} )$. One has $\beta>0$ since $\Vert \gamma \Vert$ is increasing.
    \end{itemize}
\end{definition}

To fix some notation, for a monodromy representative and a limiting path \mbox{$\gamma:[0,a]\to\mc^d$}, note that there exists by \Cref{lem:lifting} a disjoint set of paths $\wt{s_1},\dots,\wt{s_k}:[0,a]\to V$ lifting~$\gamma$ through $\wt{\pi}\vert_V$ and similarly a disjoint set of paths $p_1,\dots,p_j:[0,a]\to R$ lifting $\pi\circ\gamma$ through~$\pi\vert_R$. At any $t\in [0,a]$, denote $\{ \wt{s}_i(t) \}_{i=1}^k$ by $\wt{S}(t)$ and similarly $\{p_i(t)\}_{i=1}^j$ by $P(t)$. 

\begin{remark}\label{rem:mon_structure}
    Suppose that $\gamma$ is a limiting path for a monodromy representative.  For any \hbox{$t\in [0,a]$}, if $L_1,L_2,\dots,L_{d-1}$ are hyperplanes in $\mc^d$ defined by the vanishing of each corresponding component function of $\pi = (L_1,L_2,\dots,L_{d-1})$ and 
    $\gamma(t)\in L_i$ for $i=1,\dots,d-2$, then 
    $P(t) = R\cap \sL_{\theta(t)}$ where $\theta(t) = L_{d-1}(\gamma(t))$.  In particular, the indicated intersection of $R$ with hyperplanes is finite, and $\sL_{\theta(t)}$ is homeomorphic to a plane in 2 real dimensions. The fundamental group $\pi_1((\wt{\pi}(V)\setminus R)\cap \sL_{\theta(t)}, \gamma(t))$ is subsequently generated by homotopy classes of $j$ loops in ($\wt{\pi}(V)\setminus R)\cap \sL_{\theta(t)}$ based at $\gamma(t)$, each of which encircles exactly 1 distinct point in $P(t)$, and this fundamental group acts on the fiber $\wt{S}(t) = \wt{\pi}\vert_V^{-1}(\gamma(t))$ by monodromy.
\end{remark}

\begin{definition}\label{def:basic_loops}
    With notation as in \Cref{rem:mon_structure}, call a loop $\ell_i:[0,1]\to (\wt\pi(V)\setminus R)\cap \sL_{\theta(t)}$ a \emph{basic loop} for $p_i(t)$ if it is the concatenation of a straight line path, a path winding once counterclockwise around a circle centered at $p_i(t)$ which encircles no other point in~$P(t)$, and the reverse of the first straight line path. 
\end{definition}

\begin{remark}
    Given a localized monodromy representative for a germ $\bo{V}$ with projections $\wt{\pi}\vert_V$, $\pi\vert_R$, and limiting path $\gamma$, note that since $0\in\mc^N$ and $0\in\mc^d$ are the only elements in the fibers of $\wt{\pi}\vert_V$ and $\pi\vert_R$ over $0$ respectively, we must have that $\wt{S}(0) = \{ 0 \}$ and $P(0) = \{0 \}$.
\end{remark}

\subsection{Monodromy functors}
We will see that a monodromy representative together with a limiting path induces a functor of monodromy actions along the path. For any interval $I\subseteq \mr$, let, by abuse of notation,~$I$ also 
denote the corresponding category obtained from the poset $(I,\leq)$ where $\leq$ is the standard order. Let $\act$ be the category of group actions on sets. In the following, we suppress giving an explicit symbol for a group action where it is clear from context. More precisely, $\act$ is the category where:
\begin{itemize}
    \item If $G$ is any group and $S$ any set, any group action of $G$ on $S$ is an object of $\act$, which is denoted by $(G,S)$.
    \item An arrow $(G_1,S_1)\to(G_2,S_2)$ in $\act$ is a pair $(h,\iota)$ where $h:G_1\to G_2$ is a homomorphism and $\iota:S_1\to S_2$ is a map such that $h(g)\iota(s) = \iota(gs)$ for all $g\in G_1,s\in S_1$.
    \item Composition and identities are component-wise.
\end{itemize}
\begin{remark}
   There is a functor $\im:\act\to\group$. For any object $(G,S)$ of $\act$, view the group action as a homomorphism $\nu:G\to\text{Aut}(S)$. Then, $\im(G,S)$ is $\im(\nu)$. Given an arrow $(h,\iota):(G_1,S_1)\to(G_2,S_2)$, the homomorphism $\im(h,\iota):\im(\nu_1)\to\im(\nu_2)$ is $\nu_2\circ h\circ\nu_1^{-1}$.
\end{remark}

\begin{definition}
    Given a monodromy representative with a limiting path $\gamma$, the (unsliced) \emph{monodromy functor} for this representative, denoted by $\loc_\gamma:(0,a]\to\act$ where the monodromy representative is clear from context, is defined by the monodromy action $\loc_\gamma(t) = (\pi_1(\wt{\pi}(V)\setminus R,\gamma(t)),S(t))$. 
    For any $t_1\leq t_2 \in (0,a]$, the map $\loc_\gamma(t_1\leq t_2)$ is defined on the fiber by $s_i(t_1)\mapsto s_i(t_2)$ for $s_1(t_1),\dots,s_k(t_1)\in S(t)$ and on the fundamental group by the standard change-of-basepoint isomorphism\footnote{We adopt the typical convention here that an overline denotes the reverse of a path and $\cdot$ denotes concatenation of paths.} induced by $\ell\mapsto \overline{\gamma\vert_{[t_1,t_2]}} \cdot \ell \cdot \gamma\vert_{[t_1,t_2]}$ for all loops $\ell$ based at $\gamma(t_1)$. The \emph{sliced monodromy functor} for the monodromy representative, $\loc_\gamma^s:(0,a]\to\act$, is the same as $\loc_\gamma$ but replacing $\pi_1(\wt{\pi}(V)\setminus R,\gamma(t))$ with the image of the map on fundamental groups induced by the inclusion $(\wt{\pi}(V)\cap\sL_{\theta(t)})\setminus R \hookrightarrow \wt{\pi}(V)\setminus R$.
\end{definition}

\begin{proposition}
    $\loc_{\gamma}$ and $\loc^s_\gamma$ are functors.
\end{proposition}
\begin{proof}
We prove the unsliced case, from which the sliced case follows. The only property that is non-routine to check is that $\loc_{\gamma}(t_1\leq t_2)$ as defined is a map of group actions for any $t_1\leq t_2\in (0,a]$. Denote that map as $(\iso,\iota)$. For any $[\ell]$ in the fundamental group component of $\loc_\gamma(t_1)$, let $\Tilde{\ell}$ be the unique lifting of $\ell$ through $\wt{\pi}\vert_V$. For any $\wt s_i(t_1)\in \wt S(t_1)$, note that $\Tilde{\ell}$ ends at $[\ell]\wt 
s_i(t_1)$, which we denote by $\wt s_z(t_1)$. Also, note that $\overline{\wt s_i\vert_{[t_1,t_2]}}
\cdot\Tilde{\ell}\cdot\wt s_z\vert_{[t_1,t_2]}$ lifts $\overline{\gamma\vert_{[t_1,t_2]}}\cdot\ell\cdot\gamma\vert_{[t_1,t_2]}$, starts at $\wt s_i(t_2)$, and ends at $\wt 
s_z(t_2)$. This yields 
\[ \iota([\ell]\wt s_i(t_1)) = \iota(\wt s_z(t_1)) = \wt s_z(t_2) = \iso([\ell])\wt s_i(t_2)\text{.}\] 
\end{proof}

\begin{definition}
    The \emph{local monodromy action} of a monodromy representative with a 
    limit\-ing path 
    is $\varprojlim \loc_\gamma$ provided that the limit exists. 
    Similarly, the \emph{local monodromy group} is~$\varprojlim \im\circ\loc_\gamma$. 
\end{definition}

These limiting procedures can be understood as ``filling in'' the missing value of $\loc_\gamma$ 
at~$0$. Note that $\group$ is a complete category, so that local monodromy groups always exist. The behavior of sliced monodromy functors is even more straightforward.

\begin{figure}
    \centering
    \includegraphics[width=1.0\linewidth]{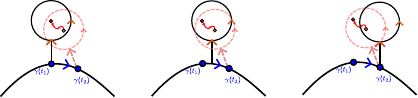}
    \begin{tabular}{ccc}
    (a) \hspace{0.95in} & \hspace{0.95in} (b) \hspace{0.95in} & \hspace{0.95in} (c)  \end{tabular}
    \caption{Illustrating first translation homotopy in \Cref{lem:generators}: (a) loop $\ell_i^{t_1}$ (black, orange arrow), $\ell_i^{t_2}$ (pink, dashed), path $p_i\vert_{[t_1,t_2]}$ (red), and path $\gamma$ (bottom, blue arrow);
    (b) loop $\ell_i^{t_1}$ has been translated to $\ell_i^{t_1}+c_t$ for some $t\in [t_1,t_2]$;
    (c) loop ends translation at $\ell_i^{t_1}+c_{\gamma(t_2)-\gamma(t_1)}$. }
    \label{fig:iso_first_translation}
\end{figure}

\begin{lemma}
For any $t_1\leq t_2$, let $\iso$ be the homomorphism component of 
\mbox{$\loc^s_\gamma(t_1\leq t_2)$}. If~$\ell_i^{t_1}$ is a basic loop for $p_i(t_1)$ and $\ell_i^{t_2}$ is a basic loop for $p_i(t_2)$, then $\iso([\ell_i^{t_1}]) = [\ell_i^{t_2}]$. In particular, $\loc^s_\gamma(t_1\leq t_2)$ is an isomorphism.    
\label{lem:generators}
\end{lemma}
\begin{proof}
    For any $t\in (0,a]$, let $d(t)$ denote the minimum of 
    both $\min_{i\not=z} \{ \Vert p_i(t) - p_z(t) \Vert\}$ and $\min \{ \Vert p_i(t) \Vert \}_{i=1}^j$, and let $D= \min_{t\in [t_1,t_2]} d(t)$. Note that $D > 0$ as it is the minimum of a continuous positive function on a compact interval. Without loss of generality, we can assume the radius of the loops $\ell_i^{t_1}$ and $\ell_i^{t_2}$ is less than $D/2$. 
    
    First, consider the simpler case where $p_i(t)\in B_{D/2}(p_i(t_1)+\gamma(t)-\gamma(t_1))$ for all $t\in [t_1,t_2]$. If necessary, shrink $D$ further such that the indicated ball does not contain $0$ or any point~$p_z(t)$ for $z\not= i$ and any $t\in [t_1,t_2]$. In this case, there is a basepoint-preserving homotopy of loops $H:I\times I \to (\wt{\pi}(V)\setminus R)$ between $\iso(\ell_i(t_1))$ and $\ell_i(t_2)$. It is given in two steps. For the first, $H_1:I\times [t_1,t_2]$ let $c_t$ be the constant loop at $\gamma(t)-\gamma(t_1)$. Then
    \[ 
    H_1(\bullet,t) = \overline{\gamma\vert_{[t,t_2]}} \cdot (\ell_i^{t_1} + c_t)\cdot\gamma\vert_{[t,t_2]} \text{.}
    \]
    See \Cref{fig:iso_first_translation}. At the end of this homotopy, $H_1(\bullet,1)$ is homotopic to the translated loop $\ell_i^{t_1}+c_{t_2}$ with circular portion of radius less than $D/2$ centered at $p_i(t_1)+\gamma(t_2)-\gamma(t_1)$, the ball $B_{D/2}(p_i(t_1)+\gamma(t_2)-\gamma(t_1))\cap \sL_{\theta(t_2)}$ contains $p_i(t_2)$, and it contains no other point $p_z(t_2)$. Note that $\ell_i^{t_1}+c_{t_2}$ and $\ell_i^{t_2}$ have images contained in $(\wt{\pi}(V)\cap\sL_{\theta(t_2)})\setminus R \cong B \setminus F$ for some closed ball $B\subseteq\mr^2$ and finite set $F$. It is subsequently straightforward to see that $\ell_i^{t_1}+c_{t_2}$ and $\ell_i^{t_2}$ are homotopic in $\wt{\pi}(V)\setminus R$.

    For the general case, notice that since $t\mapsto \Vert p_i(t) - (p_i(t') + \gamma(t)-\gamma(t_1))\Vert$ is continuous for fixed $t'\in [t_1,t_2]$, for every $t\in [t_1,t_2]$ there an open neighborhood of $t$ in $[t_1,t_2]$ for which the simple case applies. Since $[t_1,t_2]$ is compact, the first part of the lemma follows from finitely many applications of the simpler case. The map on actions given by 
    \mbox{$\wt s_i(t_2)\mapsto \wt s_i(t_1)$} for \mbox{$i=1,\dots,k$} and induced by $[\ell_i(t_2)]\mapsto [\ell_i(t_1)]$ for $i=1,\dots,j$ is therefore an inverse~to~$\loc_\gamma^s(t_1\leq t_2)$.
\end{proof}
It follows directly that the inverse limit of a sliced monodromy functor $\loc_\gamma^s$ exists and is isomorphic to $\loc_{\gamma}^s(t)$ for any $t$ in the domain of $\gamma$. We find the 
following~observations~useful.

\begin{proposition}\label{prop:slice_limits}
    If $\gamma_1:[0,a_1]\to\mc^d$ and $\gamma_2:[0,a_2]\to\mc^d$ are limiting paths for the same monodromy representative, then $\varprojlim \loc_{\gamma_1}^s \cong \varprojlim \loc_{\gamma_2}^s$.
\end{proposition}
\begin{proof}
    There exist some $t_1,t_2\in [0,\min(a_1,a_2)]$ with $\Vert \gamma_1(t_1)\Vert \leq \Vert \gamma_2(t_2)\Vert$ since the paths are increasing in norm. Form a new limiting path $\gamma_3:[0,a_3]\to\mc^d$ that ends at $\gamma_2(t_2)$ and includes $\gamma_1(t_1)$ in its image. By \Cref{lem:generators}, all the maps $\loc^s_{\gamma_i}(t'_1\leq t'_2)$ are isomorphisms for $i=1,2,3$ with $t'_1\leq t'_2$ in the appropriate interval. The limits in question are therefore isomorphic to $\loc^s_{\gamma_i}(t_i)$ for $i=1,2$, and there is an isomorphism between those 
    defined~by~$\loc^s_{\gamma_3}$.
\end{proof}

\begin{proposition}\label{prop:func_restrictions}
    Consider a monodromy representative with a limiting path restricted to $\wt{B}_{\wt{\epsilon}}$ and $B_\rho$. Let $\gamma_1:[0,a_1]\to\mc^d$ be the original limiting path and $\gamma_2 = \gamma_1\vert_{[0,a_2]}$ be the restriction of the original limiting path so obtained.      
    If $\loc^r_{\gamma_2}$ is the (sliced or unsliced) monodromy functor of the restricted monodromy representative and $\loc_{\gamma_2}$ for the unrestricted representative with path $\gamma_2$, there is a natural transformation $\loc^r_{\gamma_2}\Rightarrow \loc_{\gamma_2}$ induced by inclusion. If the local monodromy actions exist, the natural transformation induces~a~map~$\varprojlim \loc^r_{\gamma_2} \to \varprojlim \loc_{\gamma_2}$.        
\end{proposition}
\begin{proof}
Follows from a routine checking of definitions.
\end{proof}

\subsection{Local monodromy actions of holomorphic subvarieties}
We are now in position to state and prove our 
collection of main results, 
which will allow for computations with sliced monodromy functors rather than unsliced ones. 
Using the notation of \Cref{def:mon_data}, we will assume throughout this subsection that limiting paths are restricted to codomain $\cap_{i=1}^{d-2} L_i$ in order to fulfill the requirements of the local Lefschetz-Zariski theorem~in~\Cref{lem:lefschetz}.
\begin{theorem}\label{thm:slice}
    Given any localized monodromy representative for a pure $d$-dimensional germ $\bo{V}$ of a holomorphic subvariety of $\mc^N$ with limiting path $\gamma_1:[0,a_1]\to \cap_{i=1}^{d-2} L_i$, the local monodromy action of this representative exists and is isomorphic to the sliced limit $\varprojlim \loc_{\gamma_2}^s$ for some restriction $\gamma_2 = \gamma_1\vert_{[0,a_2]}$. 
\end{theorem}
\begin{proof}
    By \Cref{lem:lefschetz} and the definition of a localized monodromy representative, there exists $\theta>0$ such that, if $\vert L_{d-1}(\gamma(t)) \vert < \theta$, then $\loc^s_{\gamma_1}(t) \cong \loc_{\gamma_1}(t)$. Since $\gamma_1$ is an increasing path, $\vert L_{d-1}\circ \gamma_1 \vert$ is a continuous increasing function. It follows that there exists $a_2 > 0$ with $\vert L_{d-1}(\gamma_1(t))\vert <\theta$ for all $t\in [0,a_2]$. From \Cref{lem:generators} one immediately has that $\varprojlim \loc^s_{\gamma_2}$ exists and is isomorphic to $\loc^s_{\gamma_2}(t)$ for any $t\in [0,a_2]$. The result now follows from \Cref{prop:slice_limits}. 
\end{proof}

\begin{definition}
    The \emph{local monodromy action} of a pure $d$-dimensional holomorphic germ~$\bo{V}$ of an open subset of $\mc^N$ is the local monodromy action of any localized monodromy representative for $\bo{V}$ along any limiting path $\gamma$. Define the \emph{local monodromy group} of $\bo{V}$ similarly. 
\end{definition}

In principle, this definition depends on the generic choices of data used to construct a localized monodromy representative for $\bo{V}$, and we must check that different choices yield isomorphic local monodromy actions. The following corollary first justifies referring to \emph{the} local monodromy action and group of a localized monodromy representative.

\begin{corollary}\label{cor:unique_paths_and_restrictions}
    For any localized monodromy representative with limiting paths 
    \[\gamma_1,\gamma_2:[0,a_i]\to\cap_{i=1}^{d-2} L_i\] 
    for a pure $d$-dimensional germ $\bo{V}$ of a holomorphic subvariety of $\mc^N$:
    \begin{enumerate} 
        \item The local monodromy action is isomorphic to $\loc^s_{\gamma_1}(t)$ for any $t\in [0,a_1]$. 
        \item The local monodromy actions defined by $\gamma_1$ and $\gamma_2$ are isomorphic. 
        \item The local monodromy actions defined by any two restrictions of the same localized monodromy representative for $\bo{V}$ with limiting path $\gamma$ are isomorphic. 
    \end{enumerate}
\end{corollary}
\begin{proof}
    The first statement follows from \Cref{thm:slice} and \Cref{lem:generators}. The second follows from \Cref{thm:slice} and \Cref{prop:slice_limits}. For the third statement, note that for any two restrictions of the same localized monodromy representative, there exists a common restriction of both with monodromy action $\loc_1$. Let $\loc_2,\loc_3$ be the monodromy actions of the other two restrictions. From the first statement, we know that there is some $t_0$ such that $\varprojlim \loc_i \cong \loc^s_i(t_0)$ for $i=1,2,3$. We can further choose $t_0$ small enough that the corresponding fibers $\wt{S_i}(t_0)$ and $P_i(t_0)$ have the same number of points for $i=1,2,3$. From \Cref{prop:func_restrictions}, we have induced maps $\loc^s_1(t_0)\to\loc^s_i(t_0)$ for $i=2,3$, and one may observe that these are isomorphisms from their definition. 
\end{proof}

It is also natural to ask whether the choice of generic projections $\wt{\pi}$ and $\pi$ for the localized monodromy representative impact the action. Handling $\pi$ is straightforward, but $\wt{\pi}$ takes some care since this determines the branch locus.

\begin{proposition}\label{prop:unique_lower_projection}
    The local monodromy actions of any two localized monodromy representatives which differ only in the linear projections $\pi^1,\pi^2:\mc^d\to\mc^{d-1}$ have isomorphic local monodromy actions.
\end{proposition}
\begin{proof}
    By~\Cref{thm:slice}, there is a limiting path $\gamma$ for both monodromy representatives and 
    $t_0>0$ where, denoting the corresponding sliced monodromy actions by $\loc^s_i$ for $i=1,2$, one 
    has that the corresponding local monodromy actions are isomorphic to $\loc^s_i(t_0)$ for $i=1,2$, 
    respectively. By standard parameter homotopy continuation results, e.g., see~\cite{morgansommese1989}, given a generic path $\Gamma:
    [0,1]\to \text{Gr}(d-1,d)^{d-1}$, i.e., $\Gamma(T) = (L_1(T),\dots,L_{d-1}(T))$, and denoting 
    \[\sL_T := \cap_{i=1}^{d-2} \{ L_i(T) = 0 \} \cap \{ L_{d-1}(T) = [L_{d-1}(T)](\gamma(t_0)) \},\] 
    the set $\{ R\cap \sL_t \}_{T\in [0,1]}$ can be parameterized as a disjoint set of smooth paths $\delta_i:
    [0,1]\to\mc^d$ for $i=1,\dots,j$. If $\Gamma$ starts at the hyperplanes defining $\pi^1$ and ends 
    at those defining $\pi^2$, a similar loop-translating argument to \Cref{lem:generators} shows that 
    the mapping from \mbox{$\loc^s_1(t_0)\to\loc^s_2(t_0)$} which is the identity on the fiber and on the fundamental groups is induced by one which maps a basic loop for $\delta_i(0)$ to a basic loop for $\delta_i(1)$ for each $i=1,\dots, j$ is an isomorphism of group actions.
\end{proof}

\begin{proposition}\label{prop:unique_upper_projection}
    The local monodromy actions of any two localized monodromy representatives which differ only in the linear projections $\wt{\pi}^1,\wt{\pi}^2:\mc^N\to\mc^d$ have isomorphic local monodromy actions.
\end{proposition}
\begin{proof}
One can follow a similar argument as in \Cref{prop:unique_lower_projection}, but now
apply parameter homotopy continuation results, e.g., see~\cite{morgansommese1989}, to a parameter homotopy between 
$\wt{\pi}^1$ and $\wt{\pi}^2$ so that the local monodromy
actions along the path between them are isomorphic.  
\end{proof}

Taken together, \Cref{cor:unique_paths_and_restrictions}
and \Cref{prop:unique_lower_projection,prop:unique_upper_projection} show that the various generic choices made when forming a localized monodromy representative yield isomorphic local monodromy actions.
\begin{theorem}
    The local monodromy actions of a pure $d$-dimesional holomorphic 
    germ~$\bo{V}$ of an open subset of $\mc^N$ defined by any two sets of localized monodromy representative data are~isomorphic. 
\end{theorem}

\begin{theorem}\label{thm:inject}
    Given any monodromy representative which restricts to a localized monodromy representative with limiting path $\gamma_1:[0,a_1]\to\cap_{i=1}^{d-2} L_i$ of a pure $d$-dimensional germ~$\bo{V}$ of a holomorphic subvariety of $\mc^N$ , consider $\wt{S}^l(t)$ and $P^l(t)$ which are the subsets of fiber points and sliced branch points respectively for $t\in [0,a_1]$ with corresponding solution paths beginning at $0$. Then, for any $t\in [0,a_1]$ the local monodromy action of $\bo{V}$ is isomorphic to the sub-action of $\loc^s_{\gamma_1}(t)$ which is comprised of the subgroup generated by homotopy classes of basic loops around the points of $P^l(t)$ acting on the points of $\wt{S}^l(t)$.
\end{theorem}
\begin{proof}
    By assumption, there is some restriction of the monodromy representative which is a localized monodromy representative for $\bo{V}$, say with restricted limiting path $\gamma_2 = \gamma_1\vert_{[0,a_2]}$, sliced monodromy action $\loc^s_{\gamma_1}$ for the original representative, and (sliced) monodromy functor~$\loc^r_{\gamma_2}$ for the restriction of the monodromy representative. For any $t\in [0,a_2]$, one may observe directly by definition that the image of the arrow $\loc^r_{\gamma_2}(t) \to \loc^s_{\gamma_2}(t) = \loc^s_{\gamma_1}(t)$ as defined in \Cref{prop:func_restrictions} is the described sub-action of $\loc^s_{\gamma_1}(t)$ and that the arrow is monic. For $t \geq a_2$, compose the arrow $\loc^r_{\gamma_2}(t)\to \loc^s_{\gamma_1}(t)$ with the isomorphism $\loc^s_{\gamma_1}(a_2\leq t)$.
\end{proof}

\section{Computing local monodromy actions}\label{sec:algorithm}

The theoretical results of \Cref{sec:theory} 
yield an approach 
for computing local monodromy actions.
In particular, this theory shows that one can use analytic
continuation to extend beyond the small enough neighborhood
restriction for localizing monodromy computations.
In the following, we specialize the setup 
to the situation in~\Cref{ex:algebraic_continuation} 
for a pure $d$-dimensional algebraic complete intersection $V = \V(f)$.
For $x^*\in V$, it is straightforward to construct a monodromy representative for a projection map germ $\boldsymbol{\wt\pi}:\bo{V}\to\bo{C}^d$ that restricts to a localized monodromy representative of $\bo{V}$. The following uses the 
notation of \Cref{ex:algebraic_continuation} 
with superscripts $g$ and $l$ denoting ``global'' and ``local,'' respectively. 
\begin{algorithm}[h!]
    \scriptsize
	\SetKwInOut{input}{Input}\SetKwInOut{output}{Output}\SetKwFunction{Return}{Return}
	\input{A polynomial system $f:\bC^N\rightarrow\bC^{N-d}$ defining a reduced complete intersection $V = V(f)\subseteq\mc^N$ of dimension $d$ and a point $\starx\in V$.}
    \output{A numerical local irreducible decomposition of $V$ at $\starx$.} 
    \smallskip
     Select linear maps $\wt{\pi}:\mc^N\to\mc^d$ and $\pi:\mc^d\to\mc^{d-1}$ from open bounded sets of maps at uniform random\;
     Select $\gamma(1)\in B_1\cap_{i=1}^{d-2} L_i$ at uniform random and set $\gamma(t) = t\gamma(1)$ for $t\in[0,1]$, denote $\theta:=\pi\circ\gamma$\;
     Set $\wt{W_t}^g:=\wt\pi^{-1}(\gamma(t))\cap V$ and $W_t^g:= CG \cap (\mc^N\times\sL_{\theta(t)})$, and compute $\wt{W_1}^g$ and $W_1^g$\;
     Use homotopy continuation to compute $\wt{W_0}^g$ and $W_0^g$ by tracking the solution paths $\wt{W_t}^g$ and 
     $W_t^g$ from $t=1$ to $t=0$ starting at $\wt{W_1}^g$ and $W_1^g$, respectively\;
     Set $\wt W^l$ and $W^l$ to be the set of points in $\wt{W_1}^g$ and $W_1^g$, respectively, whose solution paths converge to $\starx$ and $(\starx,\wt\pi(\starx))$ respectively as $t\to 0$\; 
     Compute the partition $\wt{W}^l = \wt w_1 \coprod \wt w_2 \coprod \dots \coprod \wt w_k $ given by the monodromy action on $\wt W^l$ through $\wt{\pi}$ of basic loops from the second factor of $W^l$\; 
     \Return the sets $\{ f,\wt\pi, \wt w_i\}$ for $i=1,\dots,k$ \;
	\caption{\textsc{Numerical local irreducible decomposition}}
	\label{alg:decomposition}
\end{algorithm}

\Cref{lem:local} and \Cref{thm:inject} justify the following.

\begin{theorem}\label{thm:AlgCorrect}
Assuming generic choices, \Cref{alg:decomposition} 
is correct.      
\end{theorem}

Note that restricting the input of \Cref{alg:decomposition} to 
complete intersections is not strictly necessary, but
is stated this way for simplicity of presentation. 
There are standard techniques used in numerical algebraic geometry,
e.g., randomization and Bertini's theorem~\cite[\S A.9]{SommeseWamplerBook}, 
for reducing to this case. 
Note that such reduction techniques 
simply add pre- and post-processing steps
which do not change the core procedure. 

Several software packages in numerical algebraic
geometry are able to perform the computations in the key parts
of \Cref{alg:decomposition}, namely in 
lines $3$, $4$, and $6$, including
{\tt Bertini}~\cite{BHSW06},
{\tt HomotopyContinuation.jl}~\cite{HomotopyContinuationJL},
and 
{\tt NAG4M2}~\cite{NAG4M2}.
Since the computation in \Cref{alg:decomposition} is only interesting
when $x^*$ is a singular point of $V$, 
endgames (see \cite[Chap.~10]{SommeseWamplerBook} and \mbox{\cite[Chap.~3]{BertiniBook}}) can be employed in line $4$ to accurately
compute $\wt{W_0}^g$ and $W_0^g$.  
Additionally, line~$6$ may require tracking paths near
singularities and thus can be poorly conditioned.
In particular, if the critical points in $W_1^g$ cluster closely together, the corresponding basic loops encircling points from $W^l$ must necessarily pass close to the critical locus.
One remedy is to utilize adaptive precision path tracking, 
e.g., see~\cite{AMP}. Alternatively, it can often
be computationally less expensive in practice to simply try 
again with different random choices.
If desired, the computations in line $6$ can be 
certified using certified path tracking, e.g., see~\cite{RobustCertification,CertifiedNewton},
or using an {\em a posteriori} certification approach~\cite{Aposteriori}.

\section{Examples}\label{sec:examples}

We conclude with several examples of computing a
numerical local irreducible decomposition
using an implementation of \Cref{alg:decomposition} 
available at \url{https://github.com/P-Edwards/LocalMonodromy.jl}\footnote{A static version of the package together with files suitable for reproducing the examples is available at \url{https://doi.org/10.5281/zenodo.14532556}.}. 
This implementation uses {\tt HomotopyContinuation.jl}~\cite{HomotopyContinuationJL} for path tracking without certification. Computation times are reported with homotopy continuation parallelized using an Intel Core i7-920 CPU with 8 CPU threads. 
Memory requirements were less than 2GB. To fix some convenient terminology from \Cref{alg:decomposition}, call the number of points in $\wt W_1^g$ the \emph{global fiber degree} of $V$, $\vert \wt{w_i} \vert$ the \emph{local fiber degree} of each locally irreducible component $V_i$, $\vert W_1^g\vert$ the \emph{global branch degree}, and $\vert W^l\vert$ the \emph{local branch degree}.
 
\begin{example}
\label{ex:whitney_umbrella2} 

\begin{figure}[!t]
    \centering
    \includegraphics[width=0.5\linewidth]{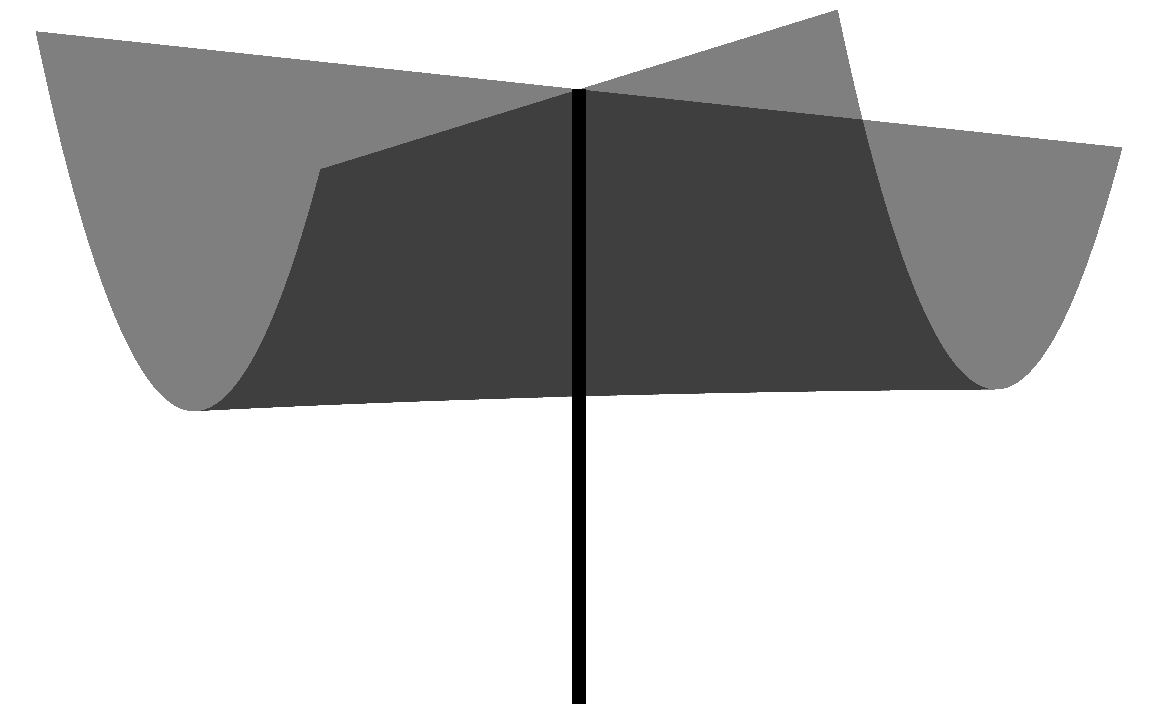}
    \caption{Illustration of the the Whitney umbrella.}
    \label{fig:Whitney}
\end{figure}    

The Whitney umbrella $V\subseteq\mc^3$ is the surface defined by $x_1^2-x_3x_2^2 = 0$ and has singular points along the line $x_1=x_2=0$
as illustrated in Figure~\ref{fig:Whitney}.
At a point $x^* = (0,0,\kappa)$ with $\kappa\not= 0$, there is a nontrivial factorization $(x_1-x_2\sqrt{x_3})(x_1+x_2\sqrt{x_3})=0$ in the local ring of holomorphic germs at $x^*$, so one expects $V$ to be locally reducible with two locally irreducible components at $x^*$. 
When $\kappa=0$, this factorization is not available since there is no holomorphic inverse to $z\mapsto z^2$ in a neighborhood of $0$. One therefore expects $V$ to be locally irreducible at the origin. 

Algorithm~\ref{alg:decomposition} at the origin
computes one local irreducible component 
of local fiber degree~2, global fiber degree 3, 
local branch degree 2, and global branch degree 4. 
The results at $x^*=(0,0,-1)$ were similar except
having two local irreducible components, each having
local fiber degree $1$.
In these cases, the computations required took approximately
22 seconds.  

This example is notable for having a critical point locus with an unreduced irreducible component, namely the line $x_1=x_2=0$ which is often called the ``handle''
of the Whitney umbrella. It is unreduced in the sense of having generic multiplicity greater
than $1$.  A standard approach in numerical algebraic geometry to perform
computations on such components is to \emph{deflate} the component
first described in~\cite[\S 10.5]{SommeseWamplerBook} (see also~\cite{Isosingular}). Our experiments include a deflation step for this and other unreduced examples.
\end{example}

\begin{example}\label{ex:hypersurfaces}
\Cref{tab:hypersurface_table} collects data 
arising from 
various globally irreducible hypersurfaces which
have an isolated singularity at the origin.  
In particular, the last one is an example of a 
so-called Brieskorn manifold~\cite{Brieskorn1966} with 
global branch degree $174$. 
\Cref{fig:brieskorn} shows how increasing degree can complicate the clustering pattern of branch points.

    \begin{table}[!h]
        \centering
        \begin{tabular}{|c|p{0.10\linewidth}|p{0.10\linewidth}|p{0.10\linewidth}|p{0.10\linewidth}|}
            \hline
            \multirow{3}{*}{Equation} & Local fiber degrees & Global fiber degree & Local branch degree & Global branch degree \\
            \hline
            $x^2+(y-1)y^2=0$ & 1,1 & 3 & 1 & 3 \\
            \hline
            $(3x+y+2z)^2x^3+(x-1)(y+z)^3=0$ & 1,2 & 5 & 1 & 10 \\
            \hline
            $xy-z^3=0$ & 2 & 3 & 2 & 4 \\
            \hline
            $x^2+y^2+z^2=0$ & 2 & 2 & 2 & 2 \\
            \hline
            $x^2+y^2+z^2+w^2=0$ & 2 & 2 & 2 & 2 \\
            \hline
            $z_1^2+z_2^2+z_3^2+z_4^3+z_5^{59}=0$ & 2 & 59 & 2 & 174 \\
            \hline
        \end{tabular}
        \caption{Summary of results for several hypersurfaces at the origin.}
        \label{tab:hypersurface_table}
    \end{table}

\begin{figure}[!t]
    \centering
    \includegraphics[width=0.4\linewidth]{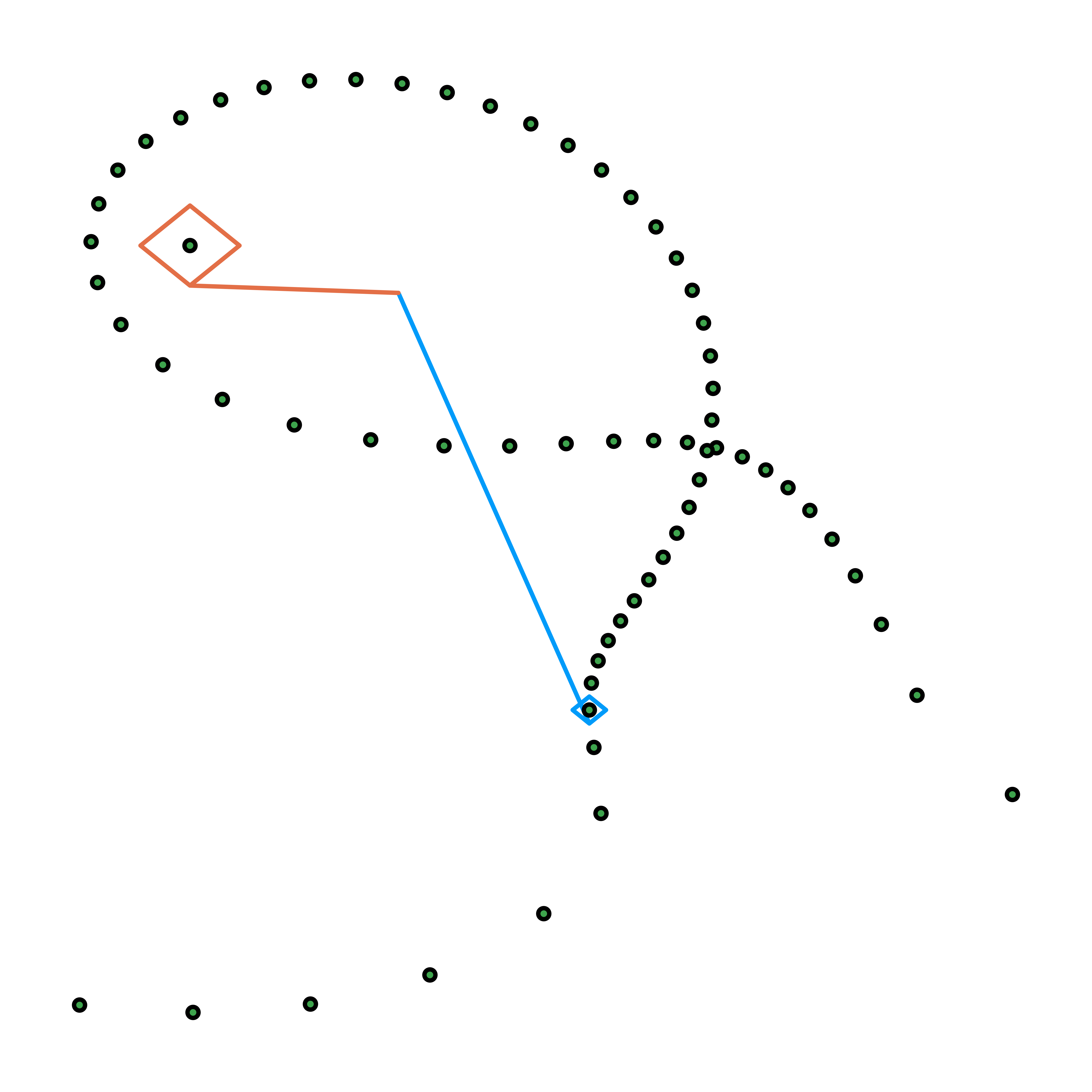}
    \caption{Branch points and corresponding localized monodromy loops of the branch locus intersected with a complex line, identified with $\mr^2$, for the Brieskorn manifold with global fiber degree 174. For illustration purposes, only 65 of 174 global branch points are depicted.}
    \label{fig:brieskorn}
\end{figure}

\end{example}

\begin{example}
Our final example arises in kinematics 
and studied in~\cite[\S8.1]{cobian2024robust}
arising from a coupler curve of a four-bar linkage.
In particular, for the polynomial system
$$f(x,y,u,v,a,b,c,d) = \left[\begin{array}{c}
    x^2+y^2-a^2 \\
    (u-b)^2+v^2-c^2 \\
    (x-u)^2+(y-v)^2-d^2 
    \end{array}\right],$$
the set $V = V(f)\subset\mc^8$ is irreducible of codimension $3$.
The variables $a,b,c,d$ are mechanical parameters of the four-bar linkage
while $x,y,u,v$ describe the coupler curve of the resulting four-bar
linkage.  We consider the local irreducible decomposition 
of $V$ at the origin.  
Here, the critical point locus has dimension $4$.
The computation found that the origin is locally irreducible
with local degree $8$ and both the global and local branch degrees 
were $24$.  Moreover, the corresponding local monodromy group was computed to be the entire symmetric group on the 8 fiber points. 
In total, this computation took approximately 45 seconds.
\end{example}

\section{Conclusion}\label{sec:conclusion}

This paper introduced a theory of local monodromy actions for germs of holomorphic subvarieties, their behavior under continuation, algorithms which leverage that theory to compute numerical local irreducible decompositions, and an open source software implementation for doing so in the algebraic case.  Several examples are used to demonstrate 
this novel theory for computing local monodromy actions
and numerical local irreducible decompositions.

\section*{Acknowledgments}

The authors thank Robin Pemantle for interesting discussions
regarding local irreducible decompositions. 
J.D.H. was supported in part by NSF CCF 2331400, Simons Foundation SFM-00005696, and the Robert and Sara Lumpkins Collegiate Professorship.

\bibliographystyle{abbrv}
\bibliography{References}

\end{document}